\documentclass[11pt,english]{article}

\usepackage[latin9]{inputenc}
\usepackage[a4paper]{geometry}
\usepackage{amsmath}
\usepackage{amssymb}
\usepackage{setspace}
\makeatletter

\providecommand{\tabularnewline}{\\}

\@ifundefined{date}{}{\date{}}
\usepackage{array}
\usepackage{mathrsfs}
\usepackage{color}

\makeatother

\usepackage{babel}
\begin{document}

\title{\textbf{A Sequence of Nested Exponential Random Variables}\\
\textbf{with Connections to Two Constants of Euler}}

\author{Michael R. Powers\thanks{Corresponding author; B306 Lihua Building, Department of Finance,
School of Economics and Management, and Schwarzman College, Tsinghua
University, Beijing, China 100084; email: powers@sem.tsinghua.edu.cn.}}

\date{March 29, 2023}
\maketitle
\begin{abstract}
\begin{singlespace}
\noindent We investigate a recursively generated sequence of random
variables that begins with an Exponential random variable with parameter
(i.e., inverse-mean) 1, and continues with additional Exponentials,
each of whose random parameter possesses the distribution of the prior
term in the sequence. Although such sequences enjoy some (limited)
applicability as models of parameter uncertainty, our present interest
is primarily theoretical. Specifically, we observe that the implied
sequence of distribution functions manifests a surprising connection
to two well-known mathematical quantities first studied by Leonhard
Euler: the Euler-Gompertz and Euler-Mascheroni constants. Through
a close analysis of one member of this distribution-function sequence,
we are able to shed new light on these constants.\medskip{}

\noindent \textbf{Keywords:} Exponential distribution; random sequence;
Euler-Gompertz constant; Euler-Mascheroni constant; irrationality.
\end{singlespace}
\end{abstract}

\section{Introduction}

\noindent We investigate the sequence of random variables, $\left\{ Y_{n}\right\} =Y_{1},Y_{2},\ldots$,
generated recursively by setting: (1) $Y_{1}\sim\textrm{Exponential}\left(\lambda_{1}=1\right)$;
and (2) $Y_{n}\sim\textrm{Exponential}\left(\lambda_{n}=Y_{n-1}\right)$
for $n\in\left\{ 2,3,\ldots\right\} $; where $\lambda_{n}\in\left(0,\infty\right)$
denotes the inverse of the Exponential mean (i.e., $Y_{n}\sim F_{Y_{n}}\left(y\right)=1-e^{-\lambda_{n}y},\:y\in\left(0,\infty\right)$).
In short,
\[
Y_{1}\sim F_{Y_{1}}\left(y\right)=1-e^{-y};
\]
\[
Y_{2}\sim F_{Y_{2}}\left(y\right)=E_{Y_{1}}\left[1-e^{-\lambda_{2}y}\mid\lambda_{2}=Y_{1}\right];
\]
\[
Y_{3}\sim F_{Y_{3}}\left(y\right)=E_{Y_{1}}\left[E_{Y_{2}\mid\lambda_{2}}\left[1-e^{-\lambda_{3}y}\mid\lambda_{3}=Y_{2}\right]\mid\lambda_{2}=Y_{1}\right];
\]
\begin{equation}
Y_{4}\sim F_{Y_{4}}\left(y\right)=E_{Y_{1}}\left[E_{Y_{2}\mid\lambda_{2}}\left[E_{Y_{3}\mid\lambda_{3}}\left[1-e^{-\lambda_{4}y}\mid\lambda_{4}=Y_{3}\right]\mid\lambda_{3}=Y_{2}\right]\mid\lambda_{2}=Y_{1}\right];\:\ldots.
\end{equation}

Such models, usually with fairly small values of $n$, can serve as
simple illustrations of parameter uncertainty. In actuarial finance,
for example, they may be used to show the impact of risk heterogeneity
on an insurance loss, $Y_{2}\mid\lambda_{2}\sim\textrm{Exponential}\left(\lambda_{2}\right)$,
generated by a single member of a continuum of exposure units with
distinctly different mean losses, $1/\lambda_{2}$. Assuming the particular
exposure unit responsible for $Y_{2}$ is selected randomly \textendash{}
and in such a way that $\lambda_{2}\sim\textrm{Exponential}$ with
parameter 1 \textendash{} then yields $Y_{2}\sim F_{Y_{2}}\left(y\right)=1-1/\left(y+1\right)$,
the cumulative distribution function (CDF) of the $\textrm{Pareto 2}\left(\alpha=1,\theta=1\right)$
distribution. Comparing the heavy-tailed unconditional random variable
$Y_{2}$ to the much lighter-tailed $Y_{2}\mid\lambda_{2}\sim\textrm{Exponential}\left(\lambda_{2}\right)$
(for fixed $\lambda_{2}$) reveals that risk-heterogeneity mixtures
can have a profound effect on the volatility of insurance losses.

In the present study, our interest in primarily theoretical. As will
become apparent in Section 2, the sequence of distribution functions,
$\left\{ F_{Y_{n}}\left(y\right)\right\} $, displays a surprising
connection to two well-known mathematical quantities first studied
by Leonhard Euler: the Euler-Gompertz constant, $\delta=0.596347\ldots$,
and the Euler-Mascheroni constant, $\gamma=0.577216\ldots$. Through
a close analysis of one member of this sequence, we are able to shed
new light on both $\delta$ and $\gamma$.

\section{Distribution Functions}

\noindent Given the nested conditional expected values presented in
(1), one can work out corresponding parameter-free expressions for
the CDFs, $F_{Y_{n}}\left(y\right)$, as follows:
\[
F_{Y_{1}}\left(y\right)=1-e^{-y};
\]
\[
F_{Y_{2}}\left(y\right)={\displaystyle \int_{0}^{\infty}}\left(1-e^{-y_{1}y}\right)e^{-y_{1}}dy_{1};
\]
\[
F_{Y_{3}}\left(y\right)={\displaystyle \int_{0}^{\infty}}\left[{\displaystyle \int_{0}^{\infty}}\left(1-e^{-y_{2}y}\right)y_{1}e^{-y_{1}y_{2}}dy_{2}\right]e^{-y_{1}}dy_{1};
\]
\[
F_{Y_{4}}\left(y\right)={\displaystyle \int_{0}^{\infty}}\left[{\displaystyle \int_{0}^{\infty}}\left[{\displaystyle \int_{0}^{\infty}}\left(1-e^{-y_{3}y}\right)y_{2}e^{-y_{2}y_{3}}dy_{3}\right]y_{1}e^{-y_{1}y_{2}}dy_{2}\right]e^{-y_{1}}dy_{1};\:\ldots.
\]

Unfortunately, simple analytic forms for $F_{Y_{n}}\left(y\right)$
rapidly become unattainable as $n$ increases, especially for odd
values of $n$. This is evident from the expression $F_{Y_{3}}\left(y\right)=-ye^{y}\textrm{Ei}\left(-y\right)$,
which must be written in terms of the exponential-integral function,
$\textrm{Ei}\left(-y\right)=-{\textstyle \int_{y}^{\infty}t^{-1}e^{-t}}dt$.
Nevertheless, the recursive derivation of $\left\{ Y_{n}\right\} $
permits a convenient reformulation. Namely, we see that
\[
F_{Y_{n}}\left(y\right)=E_{Y_{n-1}}\left[\Pr\left\{ X\leq\left(Y_{n-1}\right)y\right\} \right]
\]
\[
=\Pr\left\{ \dfrac{X}{Y_{n-1}}\leq y\right\} 
\]
for $n\geq2$, where $X\sim\textrm{Exponential}\left(\lambda=1\right)$.
Backward iteration then implies
\begin{equation}
Y_{n}\equiv\begin{cases}
{\displaystyle \prod_{j=1}^{\left(n+1\right)/2}}X_{j}\:\:/{\displaystyle \prod_{j=\left(n+3\right)/2}^{n}}X_{j} & \textrm{for }n\in\left\{ 1,3,\ldots\right\} \\
{\displaystyle \prod_{j=1}^{n/2}}X_{j}\:\:/{\displaystyle \prod_{j=\left(n/2\right)+1}^{n}}X_{j} & \textrm{for }n\in\left\{ 2,4,\ldots\right\} 
\end{cases},
\end{equation}
where the $X_{j}\sim\textrm{i.i.d. Exponential}\left(\lambda=1\right)$.

The expressions for $\left\{ Y_{n}\right\} $ in (2) reveal two important
aspects of this sequence. First, the random variable $Y_{n}$ possesses
the same probability distribution as its inverse ($1/Y_{n}$) for
all even $n$, with both distributions ``symmetric'' about $y=1$
in the sense that $F_{Y_{n}}\left(1\right)=F_{1/Y_{n}}\left(1\right)=1/2$.
Second, $Y_{n}$ is somewhat ``top-heavy'' for odd $n$, with $F_{Y_{n}}\left(1\right)>1/2$.
Computed values of $F_{Y_{n}}\left(1\right)$ provided by Table 1
show that the impact of the additional Exponential random variable
in the numerator of $Y_{n}$ decreases over $n$, with $F_{Y_{n}}\left(1\right)\rightarrow1/2$
as $n\rightarrow\infty$. The indicated values of $F_{Y_{n}}\left(1\right)$
manifest an additional, somewhat intriguing, property: both $F_{Y_{3}}\left(1\right)=0.596347\ldots$
and $F_{Y_{5}}\left(1\right)=0.577216\ldots$ are well-known mathematical
constants first studied by Euler that have long resisted attempts
to prove their (likely) irrationality.\footnote{Interestingly, it was shown by Aptekarev (2009) that at least one
of the two constants must be irrational, and by Rivoal (2012) that
at least one of them must be transcendental. }
\noindent \begin{center}
Table 1. CDF of $Y_{n}$ Evaluated at $y=1$
\par\end{center}

\noindent \begin{center}
\begin{tabular}{|c|c|}
\hline 
$n$ & $F_{Y_{n}}\left(1\right)=\Pr\left\{ Y_{n}\leq1\right\} $\tabularnewline
\hline 
\hline 
$1$ & $0.632120\ldots$\tabularnewline
\hline 
$2$ & $0.5$\tabularnewline
\hline 
$3$ & $0.596347\ldots$\tabularnewline
\hline 
$4$ & $0.5$\tabularnewline
\hline 
$5$ & $0.577215\ldots$\tabularnewline
\hline 
$6$ & $0.5$\tabularnewline
\hline 
$7$ & $0.566094\ldots$\tabularnewline
\hline 
$8$ & $0.5$\tabularnewline
\hline 
\end{tabular}\medskip{}
\par\end{center}

The first of these numbers, commonly denoted by $\delta$, is often
called the Euler-Gompertz constant (after Benjamin Gompertz) because
of its role in mortality/survival analysis \textendash{} although
Gompertz appears not to have investigated its properties. The second
number, usually denoted by $\gamma$, is often called the Euler-Mascheroni
constant (after Lorenzo Mascheroni) because Mascheroni published a
calculation of the constant that provided more decimal places than
previously known at the time. (See Lagarias, 2013.) The two constants
possess similar forms as definite integrals; that is,
\[
\delta={\displaystyle \int_{0}^{\infty}\ln\left(t+1\right)e^{-t}dt}
\]
and
\[
\gamma=-{\displaystyle \int_{0}^{\infty}\ln\left(t\right)e^{-t}dt},
\]
respectively; and are related by the equation
\begin{equation}
\delta=e\left(-\gamma+{\displaystyle \sum_{k=1}^{\infty}\dfrac{\left(-1\right)^{k+1}}{k\cdot k!}}\right),
\end{equation}
discovered by G. H. Hardy.\footnote{The infinite sum in (3) may be rewritten as $e^{-1}{\textstyle \sum_{k=0}^{\infty}}{\displaystyle \left[\left(k+1\right)^{-2}e^{1}\left(-1\right)^{k}/k!\right]}=e^{-1}E_{k\mid\mu=-1}\left[\left(k+1\right)^{-2}\right]$,
where $k$ has a quasi-Poisson distribution with mean $\mu=-1$ (characterized
by negative point masses at all $k\in\left\{ 1,3,\ldots\right\} $).
This implies $E_{k\mid\mu=-1}\left[\left(k+1\right)^{-2}\right]=\delta+\gamma e=e\left(\gamma-\textrm{Ei}\left(-1\right)\right)$,
which resembles the more conventional expression, $E_{k\mid\mu=1}\left[\left(k+1\right)^{-2}\right]=e^{-1}\left(\textrm{Ei}\left(1\right)-\gamma\right)$. } Given the observations from Table 1, we now know that the two values
also lie next to one another in an infinite sequence of (apparently
irrational) constants, $\kappa_{n}=F_{Y_{n}}\left(1\right)$, for
$n\in\left\{ 1,3,\ldots\right\} $.

As noted above, the random variable $Y_{2}$ (which can be interpreted
as the ratio of two independent $\textrm{Exponential}\left(\lambda=1\right)$
random variables), possesses a $\textrm{Pareto 2}\left(\alpha=1,\theta=1\right)$
distribution, and therefore an infinite mean. As a result, all $Y_{n}$
for $n>2$ are similarly heavy-tailed, which hinders closer examination
of their distributional properties through moment calculations and
the central limit theorem (CLT). For that reason, we now transform
$\left\{ Y_{n}\right\} $ to the logarithmic scale, working with the
sequence $\left\{ W_{n}\right\} $, where $W_{n}\equiv\ln\left(Y_{n}\right)$
for $n\in\left\{ 1,2,\ldots\right\} $.

\section{Log Transformation}

\subsection{Asymptotic Distributions}

\noindent In conjunction with sidestepping the problem of heavy tails,
transforming $Y_{n}$ to the log scale enables us to work with sums,
rather than products, of independent random variables. In particular,
we can write
\begin{equation}
W_{n}\equiv\begin{cases}
{\displaystyle \sum_{j=1}^{\left(n+1\right)/2}}\ln\left(X_{j}\right)-{\displaystyle \sum_{j=\left(n+3\right)/2}^{n}}\ln\left(X_{j}\right) & \textrm{for }n\in\left\{ 1,3,\ldots\right\} \\
{\displaystyle \sum_{j=1}^{n/2}}\ln\left(X_{i}\right)-{\displaystyle \sum_{j=\left(n/2\right)+1}^{n}}\ln\left(X_{j}\right) & \textrm{for }n\in\left\{ 2,4,\ldots\right\} 
\end{cases},
\end{equation}
where the $-\ln\left(X_{j}\right)\sim\textrm{i.i.d. Gumbel}$ with
mean $\gamma$ and variance $\pi^{2}/6$.\footnote{In the absence of standardized notation for the Gumbel distribution,
we will use $-\ln\left(X\right)\sim\textrm{Gumbel}\left(\gamma,\pi^{2}/6\right)\Longleftrightarrow X\sim\textrm{Exponential}\left(\lambda=1\right)$
hereafter.} A straightforward application of the CLT then yields the following
result.\medskip{}

\noindent \textbf{Proposition 1:} For the sequences of random variables,
$\left\{ Y_{n}\right\} $ and $\left\{ W_{n}\right\} $ (defined in
(2) and (4), respectively),
\[
\dfrac{\ln\left(Y_{n}\right)}{\pi\sqrt{n/6}}\equiv\dfrac{W_{n}}{\pi\sqrt{n/6}}\stackrel{d}{\rightarrow}\begin{cases}
\textrm{Normal}\left(-\gamma,\:1\right) & \textrm{for }n\in\left\{ 1,3,\ldots\right\} \\
\textrm{Normal}\left(0,\:1\right) & \textrm{for }n\in\left\{ 2,4,\ldots\right\} 
\end{cases}
\]
as $n\rightarrow\infty$.\medskip{}

\noindent \textbf{Proof:} For $n\in\left\{ 1,3,\ldots\right\} $,
it is clear from (4) that
\[
W_{n}\equiv-{\displaystyle \sum_{j=1}^{\left(n+1\right)/2}V_{j}}+{\displaystyle \sum_{j^{\prime}=1}^{\left(n-1\right)/2}V_{j^{\prime}}},
\]
where the $V_{j}$ and $V_{j^{\prime}}$ are i.i.d. $\textrm{Gumbel}\left(\gamma,\pi^{2}/6\right)$
random variables. Letting $W_{n}^{\left(1\right)}=-{\displaystyle {\textstyle \sum_{j=1}^{\left(n+1\right)/2}}V_{j}}$
and $W_{n}^{\left(2\right)}={\displaystyle {\textstyle \sum_{j^{\prime}=1}^{\left(n-1\right)/2}}V_{j^{\prime}}}$,
it then follows from the CLT that $W_{1}^{\left(1\right)}/\pi\sqrt{\left(n+1\right)/12}\stackrel{d}{\rightarrow}\textrm{Normal}\left(-\gamma,1\right)$
and $W_{1}^{\left(2\right)}/\pi\sqrt{\left(n-1\right)/12}\stackrel{d}{\rightarrow}\textrm{Normal}\left(0,1\right)$.
Adding these two independent random variables then yields the desired
result. For $n\in\left\{ 2,4,\ldots\right\} $,
\[
W_{n}\equiv-{\displaystyle \sum_{j=1}^{n/2}V_{j}}+{\displaystyle \sum_{j^{\prime}=1}^{n/2}V_{j^{\prime}}},
\]
and the argument proceeds in the same way. $\blacksquare$\medskip{}

The above proposition thus reveals that $Y_{n}$ is asymptotically
Lognormal with increasing accumulations of the total probability split
equally between two regions: the distant right tail, and a small neighborhood
of 0.

\subsection{Characteristic Functions and CDFs}

\noindent Another benefit of working with sums of i.i.d. random variables
is that expressions for the characteristic function are more likely
to be tractable. In the case at hand, we obtain the following proposition.\pagebreak{}

\noindent \textbf{Proposition 2:} For the sequence of random variables,
$\left\{ W_{n}\right\} $, the corresponding characteristic functions
are given by
\begin{equation}
\varphi_{W_{n}}\left(z\right)=\begin{cases}
\left(\dfrac{\pi z}{\sinh\left(\pi z\right)}\right)^{\left(n-1\right)/2}{\displaystyle \int_{0}^{\infty}e^{-t}\left[\cos\left(\ln\left(t\right)z\right)+i\sin\left(\ln\left(t\right)z\right)\right]dt} & \textrm{for }n\in\left\{ 1,3,\ldots\right\} \\
\left(\dfrac{\pi z}{\sinh\left(\pi z\right)}\right)^{n/2} & \textrm{for }n\in\left\{ 2,4,\ldots\right\} 
\end{cases}.
\end{equation}
\medskip{}
\textbf{Proof:} From (4), it is easy to see that
\[
\varphi_{W_{n}}\left(z\right)=\begin{cases}
\left[\varphi_{-\ln\left(X\right)}\left(-z\right)\right]^{\left(n+1\right)/2}\left[\varphi_{-\ln\left(X\right)}\left(z\right)\right]^{\left(n-1\right)/2} & \textrm{for }n\in\left\{ 1,3,\ldots\right\} \\
\left[\varphi_{-\ln\left(X\right)}\left(-z\right)\right]^{n/2}\left[\varphi_{-\ln\left(X\right)}\left(z\right)\right]^{n/2} & \textrm{for }n\in\left\{ 2,4,\ldots\right\} 
\end{cases},
\]
where $-\ln\left(X\right)\sim\textrm{Gumbel}\left(\gamma,\pi^{2}/6\right)$.
Since $\varphi_{-\ln\left(X\right)}\left(-z\right)=E_{X}\left[e^{iz\ln\left(X\right)}\right]=E_{X}\left[X^{iz}\right]=\Gamma\left(1+iz\right)$,
it then follows that
\[
\varphi_{W_{n}}\left(z\right)=\begin{cases}
\left[\Gamma\left(1+iz\right)\right]^{\left(n+1\right)/2}\left[\Gamma\left(1-iz\right)\right]^{\left(n-1\right)/2} & \textrm{for }n\in\left\{ 1,3,\ldots\right\} \\
\left[\Gamma\left(1+iz\right)\right]^{n/2}\left[\Gamma\left(1-iz\right)\right]^{n/2} & \textrm{for }n\in\left\{ 2,4,\ldots\right\} 
\end{cases}
\]
\[
=\begin{cases}
\left(\dfrac{\pi z}{\sinh\left(\pi z\right)}\right)^{\left(n-1\right)/2}\Gamma\left(1+iz\right) & \textrm{for }n\in\left\{ 1,3,\ldots\right\} \\
\left(\dfrac{\pi z}{\sinh\left(\pi z\right)}\right)^{n/2} & \textrm{for }n\in\left\{ 2,4,\ldots\right\} 
\end{cases}
\]
\[
=\begin{cases}
\left(\dfrac{\pi z}{\sinh\left(\pi z\right)}\right)^{\left(n-1\right)/2}{\displaystyle \int_{0}^{\infty}e^{-t}t^{iz}dt} & \textrm{for }n\in\left\{ 1,3,\ldots\right\} \\
\left(\dfrac{\pi z}{\sinh\left(\pi z\right)}\right)^{n/2} & \textrm{for }n\in\left\{ 2,4,\ldots\right\} 
\end{cases},
\]
which are equivalent to the expressions in (5). $\blacksquare$\medskip{}

Rather fortuitously, the characteristic functions of Proposition 2
are relatively easy to invert, as shown in the next result.\medskip{}

\noindent \textbf{Proposition 3:} For the sequence of random variables,
$\left\{ W_{n}\right\} $, the corresponding CDFs are given by
\begin{equation}
F_{W_{n}}\left(w\right)=\begin{cases}
\dfrac{1}{2}+{\displaystyle \int_{0}^{\infty}{\displaystyle \dfrac{\left(\pi z\right)^{\left(n-3\right)/2}\int_{0}^{\infty}e^{-t}\sin\left(\left[w-\ln\left(t\right)\right]z\right)dt}{\left[\sinh\left(\pi z\right)\right]^{\left(n-1\right)/2}}dz}} & \textrm{for }n\in\left\{ 1,3,\ldots\right\} \\
\dfrac{1}{2}+{\displaystyle \int_{0}^{\infty}\dfrac{\left(\pi z\right)^{\left(n-3\right)/2}{\displaystyle \sin\left(wz\right)}}{\left[\sinh\left(\pi z\right)\right]^{\left(n-1\right)/2}}dz} & \textrm{for }n\in\left\{ 2,4,\ldots\right\} 
\end{cases},
\end{equation}
for $w\in\Re$.\medskip{}

\noindent \textbf{Proof:} For the case of $n\in\left\{ 1,3,\ldots\right\} $,
we apply the formula of Gil-Pelaez (1951) to (5), obtaining
\[
F_{W_{n}}\left(w\right)=\dfrac{1}{2}+\dfrac{1}{2\pi}{\displaystyle \int_{0}^{\infty}\dfrac{\left[e^{izw}\varphi_{W_{n}}\left(-z\right)-e^{-izw}\varphi_{W_{n}}\left(z\right)\right]}{iz}dz}
\]
\[
=\dfrac{1}{2}+\dfrac{1}{2\pi}{\displaystyle \int_{0}^{\infty}\left(\dfrac{\pi z}{\sinh\left(\pi z\right)}\right)^{\left(n-1\right)/2}\left\{ \dfrac{\left[\cos\left(wz\right)+i\sin\left(wz\right)\right]{\displaystyle {\textstyle \int_{0}^{\infty}}e^{-t}\left[\cos\left(\ln\left(t\right)z\right)-i\sin\left(\ln\left(t\right)z\right)\right]dt}}{iz}\right.}
\]
\[
\left.-\dfrac{\left[\cos\left(wz\right)-i\sin\left(wz\right)\right]{\displaystyle {\textstyle \int_{0}^{\infty}}e^{-t}\left[\cos\left(\ln\left(t\right)z\right)+i\sin\left(\ln\left(t\right)z\right)\right]dt}}{iz}\right\} dz
\]
\[
=\dfrac{1}{2}+\dfrac{1}{2\pi}{\displaystyle \int_{0}^{\infty}\left(\dfrac{\pi z}{\sinh\left(\pi z\right)}\right)^{\left(n-1\right)/2}\left\{ \dfrac{{\displaystyle {\textstyle \int_{0}^{\infty}}e^{-t}\left[\cos\left(\left[w-\ln\left(t\right)\right]z\right)+i\sin\left(\left[w-\ln\left(t\right)\right]z\right)\right]dt}}{iz}\right.}
\]
\[
\left.-\dfrac{{\displaystyle {\textstyle \int_{0}^{\infty}}e^{-t}\left[\cos\left(\left[w-\ln\left(t\right)\right]z\right)-i\sin\left(\left[w-\ln\left(t\right)\right]z\right)\right]dt}}{iz}\right\} dz
\]
\[
=\dfrac{1}{2}+\dfrac{1}{\pi}{\displaystyle \int_{0}^{\infty}\left(\dfrac{\pi z}{\sinh\left(\pi z\right)}\right)^{\left(n-1\right)/2}\dfrac{{\displaystyle {\textstyle \int_{0}^{\infty}}e^{-t}\sin\left(\left[w-\ln\left(t\right)\right]z\right)dt}}{z}}dz,
\]
 which may be rearranged to give the relevant expression in (6). The
much simpler derivation for $n\in\left\{ 2,4,\ldots\right\} $ proceeds
analogously. $\blacksquare$

\section{The Case of $\boldsymbol{n=3}$}

\subsection{The CDF and Related Functions}

\noindent Although the mathematical forms in (6) offer limited analytical
tractability, the case of $n=3$ warrants closer examination because
of the close connection between the relevant CDF and the Euler-Gompertz
constant. Inserting $n=3$ into (6) allows us to write
\[
F_{W_{3}}\left(w\right)=\dfrac{1}{2}+{\displaystyle \int_{0}^{\infty}{\displaystyle \dfrac{\int_{0}^{\infty}e^{-t}\sin\left(\left[w-\ln\left(t\right)\right]z\right)dt}{\sinh\left(\pi z\right)}dz}}
\]
\[
=\dfrac{1}{2}+{\displaystyle \int_{0}^{\infty}e^{-t}\int_{0}^{\infty}{\displaystyle \dfrac{\left[\sin\left(wz\right)\cos\left(\ln\left(t\right)z\right)-\cos\left(wz\right)\sin\left(\ln\left(t\right)z\right)\right]}{\sinh\left(\pi z\right)}dz}dt}
\]
\[
=\dfrac{1}{2}+\dfrac{1}{2}{\displaystyle \int_{0}^{\infty}e^{-t}\left[\dfrac{\sinh\left(w\right)-\sinh\left(\ln\left(t\right)\right)}{\cosh\left(w\right)+\cosh\left(\ln\left(t\right)\right)}\right]dt}
\]
\[
=\dfrac{1}{2}+\dfrac{1}{2}{\displaystyle \int_{0}^{\infty}e^{-t}\left(\dfrac{e^{w}-e^{-w}-t+t^{-1}}{e^{w}+e^{-w}+t+t^{-1}}\right)dt}
\]
\[
={\displaystyle \int_{0}^{\infty}\dfrac{e^{w-t}}{t+e^{w}}dt}
\]
\[
=-e^{\left(e^{w}+w\right)}\textrm{Ei}\left(-e^{w}\right),
\]
for $w\in\Re$, where the final expression is equivalent (under transformation)
to the CDF of $Y_{3}$ mentioned in Section 1.

To explore the behavior of $F_{W_{3}}\left(w\right)$ in the neighborhood
of $w=0$ (which corresponds to $F_{Y_{3}}\left(y\right)$ in the
neighborhood of $y=1$), let $\mathcal{F}_{W_{3}}\left(w\right)={\displaystyle {\textstyle \int_{-\infty}^{w}}F_{W_{3}}\left(s\right)ds}=\gamma+w-e^{\left(e^{w}\right)}\textrm{Ei}\left(-e^{w}\right)$
and note that
\[
\mathcal{F}_{W_{3}}\left(0\right)=\delta+\gamma,
\]
\[
F_{W_{3}}\left(0\right)=\delta,
\]
\[
f_{W_{3}}\left(0\right)=2\delta-1,
\]
\[
f_{W_{3}}^{\prime}\left(0\right)=5\delta-3,\:\ldots.
\]
This reveals a clear pattern of integers in the linear functions of
$\delta$ associated with the various derivatives of $\mathcal{F}_{W_{3}}\left(w\right)$
at 0, which may be summarized as
\[
\mathcal{F}_{W_{3}}^{\left(0\right)}\left(0\right)=\delta B_{0}-A_{0}+\gamma,
\]
\[
\mathcal{F}_{W_{3}}^{\left(1\right)}\left(0\right)=\delta B_{1}-A_{1}+1,\textrm{ and}
\]
\[
\mathcal{F}_{W_{3}}^{\left(k\right)}\left(0\right)=\delta B_{k}-A_{k},
\]
for $k\in\left\{ 2,3,\ldots\right\} $, where $\mathcal{F}_{W_{3}}^{\left(k\right)}\left(0\right)$
denotes the $k^{\textrm{th}}$ derivative of $\mathcal{F}_{W_{3}}\left(w\right)$
at 0, and $B_{k}$ and $A_{k}$ denote the $k^{\textrm{th}}$ Bell
and ``Gould'' numbers,\footnote{More than one sequence of integers is commonly described as the ``Gould
numbers''. Therefore, the reader is referred to Gould and Quaintance
(2007) for clarification.} respectively, for $k\in\left\{ 0,1,\ldots\right\} $.

\subsection{Taylor Series and Implications}

\noindent The preceding analysis permits construction of the Taylor-series
expansion of $\mathcal{F}_{W_{3}}\left(w\right)$ around $w=0$ as
follows:
\[
T_{W_{3}}\left(w\right)=\gamma+w+{\displaystyle \sum_{\ell=0}^{\infty}}\dfrac{\left(\delta B_{\ell}-A_{\ell}\right)w^{\ell}}{\ell!}.
\]
Replacing $\mathcal{F}_{W_{3}}\left(w\right)$ with $\overline{\mathcal{F}}_{W_{3}}\left(w\right)={\displaystyle {\textstyle \int_{w}^{\infty}}\left[1-F_{W_{3}}\left(s\right)\right]ds}=-e^{\left(e^{w}\right)}\textrm{Ei}\left(-e^{w}\right)=\mathcal{F}_{W_{3}}\left(w\right)-\gamma-w$
yields
\[
\mathcal{\overline{F}}_{W_{3}}^{\left(k\right)}\left(0\right)=\delta B_{k}-A_{k}
\]
for $k\geq0$, and thus
\begin{equation}
\overline{T}_{W_{3}}^{\left(k\right)}\left(w\right)={\displaystyle \sum_{\ell=0}^{\infty}}\dfrac{\left(\delta B_{\ell+k}-A_{\ell+k}\right)w^{\ell}}{\ell!},
\end{equation}
a somewhat simpler expression.

Before proceeding further, it is useful to investigate the asymptotic
properties of this series.\medskip{}

\noindent \textbf{Proposition 4:} The Taylor series $\overline{T}_{W_{3}}^{\left(k\right)}\left(w\right)$
converges pointwise to $\mathcal{\overline{F}}_{W_{3}}^{\left(k\right)}\left(w\right)$
for all $k\in\left\{ 0,1,\ldots\right\} $ and $w\in\Re$.\medskip{}

\noindent \textbf{Proof:} Consider the remainder term associated with
the $m^{\textrm{th}}$ partial sum of $\overline{T}_{W_{3}}^{\left(k\right)}\left(w\right)$:
\[
\rho_{m+1}\left(w\right)=c_{0}\dfrac{\left(\delta B_{m+k+1}-A_{m+k+1}\right)w^{m+1}}{\left(m+1\right)!},
\]
for some $c_{0}\leq\left|w\right|$. Asakly, et al. (2014) demonstrated
that
\begin{equation}
\delta B_{\ell}-A_{\ell}=B_{\ell}\:O\left(e^{-c_{1}\ell/\left(\ln\left(\ell\right)\right)^{2}}\right)
\end{equation}
for some $c_{1}\in\left(0,\infty\right)$, and used this fact to prove
$\underset{\ell\rightarrow\infty}{\lim}\left(A_{\ell}/B_{\ell}\right)=\delta$.
Moreover, Berend and Tassa (2010) showed that
\[
B_{\ell}<\left(\dfrac{0.792\ell}{\ln\left(\ell+1\right)}\right)^{\ell}
\]
for all $\ell>0$. Combining these results with Stirling's approximation
($\ell!\sim\sqrt{2\pi}e^{-\ell}\ell^{\ell+1/2}$) then yields
\[
\rho_{m+1}\left(w\right)=O\left(\left[\dfrac{0.792\left(m+k+1\right)e^{-c_{1}/\left(\ln\left(m+k+1\right)\right)^{2}}}{\ln\left(m+k+2\right)}\right]^{m+k+1}\left(\dfrac{ew}{m+1}\right)^{m+1}\left(m+1\right)^{-1/2}\right)
\]
\[
=O\left(\dfrac{\left[0.792e^{-c_{1}/\left(\ln\left(m+k+1\right)\right)^{2}}\right]^{m+k+1}}{\left(m+1\right)^{1/2}}\dfrac{\left(m+k+1\right)^{k}\left(ew\right)^{m+1}}{\left(\ln\left(m+k+2\right)\right)^{m+k+1}}\left(\dfrac{m+k+1}{m+1}\right)^{m+1}\right)
\]
\begin{equation}
=O\left(\dfrac{\left[0.792e^{1-c_{1}/\left(\ln\left(m+k+1\right)\right)^{2}}\right]^{m+k+1}}{\left(m+1\right)^{1/2}}\dfrac{\left(m+k+1\right)^{k}w^{m+1}}{\left(\ln\left(m+k+2\right)\right)^{m+k+1}}\right),
\end{equation}
which converges to 0 as $m\rightarrow\infty$ for all $k$ and $w$.
$\blacksquare$\medskip{}

It is clear from (9) that (7) does not converge uniformly over either
$w\in\Re$ or $k\in\left\{ 0,1,\ldots\right\} $ because the right-hand
side of (9) diverges to $\infty$ as both $w\rightarrow\infty$ and
$k\rightarrow\infty$. Moreover, it is easy to see that
\[
\underset{w\rightarrow\infty}{\lim}\:\underset{m\rightarrow\infty}{\lim}\:{\displaystyle \sum_{\ell=0}^{m}}\dfrac{\left(\delta B_{\ell+k}-A_{\ell+k}\right)w^{\ell}}{\ell!}=\underset{w\rightarrow\infty}{\lim}\mathcal{\overline{F}}_{W_{3}}^{\left(k\right)}\left(w\right)
\]
\[
=0
\]
does not imply
\[
\underset{m\rightarrow\infty}{\lim}\:\underset{w\rightarrow\infty}{\lim}\:{\displaystyle \sum_{\ell=0}^{m}}\dfrac{\left(\delta B_{\ell+k}-A_{\ell+k}\right)w^{\ell}}{\ell!}=\underset{w\rightarrow\infty}{\lim}\mathcal{\overline{F}}_{W_{3}}^{\left(k\right)}\left(w\right)
\]
\[
=0.
\]
This is because ${\textstyle \sum_{\ell=0}^{m}}\left(\delta B_{\ell+k}-A_{\ell+k}\right)w^{\ell}/\ell!$
alternates irregularly between finite positive and negative quantities,
and is unbounded with respect to $w$ as $m$ increases. Consequently,
taking limits of the expression as $w\rightarrow\infty$ results in
positive or negative $\infty$, depending on the specific value of
$m$ involved.

The next result describes the behavior of (7) as $w\rightarrow-\infty$.\medskip{}

\noindent \textbf{Proposition 5:}
\[
\underset{w\rightarrow-\infty}{\lim}\left[\left(\underset{m\rightarrow\infty}{\lim}\:{\displaystyle \sum_{\ell=0}^{m}}\dfrac{\left(\delta B_{\ell}-A_{\ell}\right)w^{\ell}}{\ell!}\right)+w\right]=\underset{w\rightarrow-\infty}{\lim}\left(\mathcal{\overline{F}}_{W_{3}}^{\left(0\right)}\left(w\right)+w\right)
\]
\[
=-\gamma.
\]
\medskip{}
\textbf{Proof:} From Proposition 4, we know that
\[
\underset{m\rightarrow\infty}{\lim}\:{\displaystyle \sum_{\ell=0}^{m}}\dfrac{\left(\delta B_{\ell}-A_{\ell}\right)w^{\ell}}{\ell!}=\mathcal{\overline{F}}_{W_{3}}^{\left(0\right)}\left(w\right)
\]
for all $w\in\Re$. Since
\[
\mathcal{\overline{F}}_{W_{3}}^{\left(0\right)}\left(w\right)=-e^{\left(e^{w}\right)}\textrm{Ei}\left(-e^{w}\right)
\]
\[
=-e^{\left(e^{w}\right)}\left(\gamma+w-{\displaystyle \sum_{k=1}^{\infty}\dfrac{\left(-1\right)^{k+1}e^{kw}}{k\cdot k!}}\right),
\]
it follows that
\[
\underset{w\rightarrow-\infty}{\lim}\left(\mathcal{\overline{F}}_{W_{3}}^{\left(0\right)}\left(w\right)+w\right)=\underset{w\rightarrow\infty}{\lim}\left(\mathcal{\overline{F}}_{W_{3}}^{\left(0\right)}\left(-w\right)-w\right)
\]
\[
=\underset{w\rightarrow\infty}{\lim}\left[-e^{\left(e^{-w}\right)}\left(\gamma-w-{\displaystyle \sum_{k=1}^{\infty}\dfrac{\left(-1\right)^{k+1}e^{-kw}}{k\cdot k!}}\right)-w\right]
\]
\[
=\underset{w\rightarrow\infty}{\lim}\left[\left(-\gamma+w+{\displaystyle \sum_{k=1}^{\infty}\dfrac{\left(-1\right)^{k+1}e^{-kw}}{k\cdot k!}}\right)-w\right]
\]
\[
=-\gamma+\underset{w\rightarrow\infty}{\lim}\:{\displaystyle \sum_{k=1}^{\infty}\dfrac{\left(-1\right)^{k+1}e^{-kw}}{k\cdot k!}}
\]
\[
=-\gamma.\:\blacksquare
\]

\section{Discussion and Further Research}

\noindent In the present study, we found that the sequence of distribution
functions, $\left\{ F_{Y_{n}}\left(y\right)\right\} $, possesses
a surprising connection to both the Euler-Gompertz constant, $\delta$,
and the Euler-Mascheroni constant, $\gamma$. Specifically, $F_{Y_{3}}\left(1\right)=\delta$
and $F_{Y_{5}}\left(1\right)=\gamma$. Using the transformation $W_{n}=\ln\left(Y_{n}\right)$,
we investigated $F_{Y_{3}}\left(y\right)$ through its counterpart,
$F_{W_{3}}\left(w\right)$ (for which $F_{W_{3}}\left(0\right)=F_{Y_{3}}\left(1\right)$).

Given the connection of $\gamma$ to $F_{Y_{5}}\left(y\right)$, it
is reasonable to believe a similar analysis of that CDF could yield
further interesting results, especially with regard to the constants
$\kappa_{n}=F_{Y_{n}}\left(1\right)$, for $n\in\left\{ 7,9,\ldots\right\} $.
Although the expression for $F_{W_{5}}\left(w\right)$ given by (6)
appears insufficiently tractable, it may be possible to approach the
problem from the opposite direction. Specifically, certain integer
sequences providing rational approximations to $\gamma$ (as the Bell
and ``Gould'' numbers do for $\delta$) might be used to construct
the power-series representation of a function, $\mathcal{\overline{F}}_{W_{5}}^{\left(k\right)}\left(w\right)$,
corresponding to $\mathcal{\overline{F}}_{W_{3}}^{\left(k\right)}\left(w\right)$.
Analytical characteristics of $\mathcal{\overline{F}}_{W_{5}}^{\left(k\right)}\left(w\right)$
then could be explored and exploited in a manner similar to Propositions
4-5.

\end{document}